\documentclass[11pt]{article}
\usepackage[a4paper,margin=27mm]{geometry}
\usepackage{amsmath,amssymb,amsthm,mathtools,booktabs,microtype}
\usepackage[T1]{fontenc}
\usepackage{lmodern}
\usepackage{xcolor}
\usepackage[colorlinks=true,linkcolor=blue!50!black,citecolor=blue!50!black,urlcolor=blue!50!black]{hyperref}
\newtheorem{theorem}{Theorem}[section]
\newtheorem{corollary}[theorem]{Corollary}
\newtheorem{lemma}[theorem]{Lemma}
\newcommand{\NN}{\mathbb N}
\newcommand{\ZZ}{\mathbb Z}
\newcommand{\CT}{\operatorname{CT}}
\newcommand{\Tr}{\operatorname{Tr}}
\newcommand{\supp}{\operatorname{supp}}
\newcommand{\vac}{\boldsymbol{0}}
\newcommand{\bA}{\mathsf A}
\newcommand{\bB}{\mathsf B}
\newcommand{\bT}{\mathsf T}
\title{Support enumerators of two-row grids and comb posets}
\author{Yueming Zhong\\[6pt]
\small Key Laboratory of Engineering Modeling and Statistical Computation of Hainan Province\\
\small School of Mathematics and Statistics, Hainan University, Haikou 570228, P.R. China\\[3pt]
\small Corresponding author: \href{mailto:zhongyueming107@gmail.com}{\texttt{zhongyueming107@gmail.com}}}
\date{}
\hypersetup{pdfauthor={Yueming Zhong},pdftitle={Support enumerators of two-row grids and comb posets}}
\begin{document}
\maketitle
\begin{abstract}
We enumerate lattice points of preorder polytopes by support size for two-row Ferrers posets and recursively defined comb posets. For two-row posets, the order-ideal inequalities give a reflected-walk model in the nonnegative quadrant. Its transfer operators have a rank-one commutator, from which we derive an identity with five independent parameters. A matching-parameter specialization factors into the enumerators of a Motzkin excursion and a Motzkin bridge and proves the conjectured product of type~A and type~B Narayana polynomials for a two-row grid with its maximum removed. The general identity also records the two row supports separately. For combs, maximum ideal excess leads to a scalar recursion and a cubic generating function with independent spine and tooth weights. This proves a compositional-inverse conjecture of Athanasiadis and Chapoton and gives explicit $\gamma$-coefficients by Lagrange inversion. We prove $\gamma$-positivity for all two-row Ferrers cell posets and for the combs. The truncated-grid enumerators have simple negative zeros.
\end{abstract}

\section{Introduction}

Throughout, $\NN=\{0,1,2,\ldots\}$. We set $\binom Nj=0$ for integers $N\ge0$ and $j\notin\{0,\ldots,N\}$. For a finite poset $P$, its preorder polytope is
\begin{equation}\label{eq:polytope}
 \mathcal{Q}_P=\left\{\mathbf x\in\mathbb R_{\ge0}^{P}:
 \sum_{e\in I}x_e\le |I|\text{ for every order ideal }I\subseteq P\right\}.
\end{equation}
Following Athanasiadis and Chapoton~\cite{AC}, write
\begin{equation}\label{eq:Hdef}
 H_P(z)=\sum_{\mathbf x\in \mathcal{Q}_P\cap\ZZ^{P}}z^{|\supp(\mathbf x)|}.
\end{equation}
Here $\supp(\mathbf x)=\{e\in P:x_e\ne0\}$. Athanasiadis and Chapoton call $H_P$ the $h$-polynomial of $P$. We use the term support enumerator to distinguish it from the Ehrhart $h^*$-polynomial of $\mathcal{Q}_P$.

Athanasiadis~\cite[Theorem 1.2]{Ath} proved $\gamma$-positivity for the support enumerators of composition polytopes, including the chain case considered below. Dai, Hou, Liu, Thawinrak, and Wang~\cite[Theorems 1.2 and 1.3, Section 4.3]{DHLTW} realized preorder support enumerators as $h^*$-polynomials of associated root polytopes and proved their invariance under duality.

Ding and Ding~\cite[Theorems 1.1 and 3.4]{DD} proved real stability for a family of bivariate antichain polynomials on products of a chain with a two-row Ferrers poset, and obtained a common interleaver for two boundary sums. Luo~\cite[Theorem 1.1]{Luo} proved multiplicativity under ordinal sums for the Ehrhart $h^*$-polynomials of polar preorder polytopes. These antichain and polar Ehrhart enumerators are in general different from $H_P$.

We give explicit evaluations and $\gamma$-expansions for two families from~\cite[Sections 9.3 and 9.6]{AC}: truncated two-row grids and recursively defined comb posets.

For two-row Ferrers posets, prefix sums and a running minimum encode the order-ideal inequalities as a reflected walk in the nonnegative quadrant. The transfer operators have a rank-one commutator. Theorem~\ref{thm:transfer} evaluates their return moments with independent diagonal weights, down-step weights, and a boundary coupling. The specialization in~\eqref{eq:transfer-product} factors into the enumerators of a Motzkin excursion and a Motzkin bridge. Substituting the support weights then gives the Narayana product. The general identity gives the separate row-support refinement.

The comb calculation uses the maximum excess of a coordinate sum over the size of an order ideal. This statistic is additive under disjoint unions and has a simple update when a new minimum is adjoined (Lemma~\ref{lem:excess}). Applied to the comb recursion, these identities give a scalar recurrence with a slack parameter. A catalytic variable records the slack, and the resulting kernel equation yields a cubic with independent spine and tooth weights.

Let $C_n$ denote the chain with $n$ elements. For $n\ge1$, set
\[
 P_n=(C_n\times C_2)\setminus\{\max\},\qquad G_n(z)=H_{P_n}(z).
\]
Define
\begin{equation}\label{eq:NandB}
 N_n(z)=\frac1n\sum_{j=0}^{n-1}\binom nj\binom n{j+1}z^j,
 \qquad B_n(z)=\sum_{j=0}^{n}\binom nj^2z^j.
\end{equation}
These are $\operatorname{Nar}^{\mathbb A}_{n-1}$ and $\operatorname{Nar}^{\mathbb B}_n$, respectively, in the notation of~\cite{AC}. Their product appears in the conjecture at the end of~\cite[Section 9.3]{AC}.

\begin{theorem}\label{thm:product}
For every $n\ge1$,
\begin{equation}\label{eq:product}
 G_n(z)=N_n(z)B_n(z).
\end{equation}
Consequently,
\begin{equation}\label{eq:total}
 |\mathcal{Q}_{P_n}\cap\ZZ^{2n-1}|=\frac1{n+1}\binom{2n}{n}^{\!2}.
\end{equation}
\end{theorem}

For integers $m\ge k\ge0$, let $P_{m,k}$ be the poset on
\[
 e_1<\cdots<e_m,\qquad f_1<\cdots<f_k,
\]
with additional relations $e_i<f_j$ whenever $i\le j$. Equivalently, $P_{m,k}$ is the Ferrers cell poset of shape $(m,k)$, defined here as the componentwise order on
\[
 \{(i,1):1\le i\le m\}\cup\{(j,2):1\le j\le k\}.
\]
In particular, $P_n=P_{n,n-1}$. We take $P_{0,0}$ to be empty and set $H_{P_{0,0}}=1$.

Recall that a polynomial of degree $D$ is $\gamma$-positive if it can be written as
\begin{equation}\label{eq:gamma-def}
 \sum_{j=0}^{\lfloor D/2\rfloor}\gamma_j z^j(1+z)^{D-2j},
 \qquad \gamma_j\ge0.
\end{equation}
Such a polynomial is palindromic and unimodal.

\begin{theorem}\label{thm:ferrers-gamma}
For all $m\ge k\ge0$, the polynomial $H_{P_{m,k}}(z)$ is $\gamma$-positive. More precisely, its $\gamma$-coefficients are nonnegative integers counting the return walks described in Section~\ref{sec:consequences}.
\end{theorem}

The second family consists of the comb posets $T_n$. Start with $T_0=\varnothing$. To obtain $T_{n+1}$, take the disjoint union of $T_n$ with a singleton and then adjoin a new minimum element below every element of that union. Thus $|T_n|=2n$. Put
\[
 K_n(z)=H_{T_n}(z),\qquad
 F(x,z)=\sum_{n\ge0}z^{-n}K_n(z)x^{n+1}.
\]
The compositional-inverse conjecture in~\cite[Section 9.6]{AC} is equivalent to the cubic equation below.

\begin{theorem}\label{thm:comb}
The series $F(x,z)$ is the unique solution in $x\mathbb Q[z,z^{-1}][[x]]$ of
\begin{equation}\label{eq:comb-cubic}
 F=x(1+F)\bigl(1+(1+z)F\bigr)\bigl(1+(1+z^{-1})F\bigr).
\end{equation}
For every $n\ge0$,
\begin{equation}\label{eq:comb-gamma-main}
 K_n(z)=\sum_{j=0}^{n}\eta_{n,j}z^j(1+z)^{2n-2j},
 \qquad
 \eta_{n,j}=\frac1{n+1}\binom{n+1}{j+1}\binom{2n+1-j}{j}.
\end{equation}
In particular, $K_n(z)$ is $\gamma$-positive, palindromic, and unimodal.
\end{theorem}

Theorems~\ref{thm:bivariate} and~\ref{thm:comb-bivariate} give refinements by the two row supports and by the spine and tooth supports.

We construct the grid walk model in Section~\ref{sec:walks} and prove the transfer identity in Section~\ref{sec:trace}. Its application to the Narayana product and the row refinement occupies Section~\ref{sec:refinement}. Section~\ref{sec:combs} treats combs with separate spine and tooth weights. The positivity and zero results are proved in Section~\ref{sec:consequences}.

\section{Boundary statistics and reflected walks}\label{sec:walks}

For a finite poset $P$ and $\mathbf a\in\NN^P$, define its maximum ideal excess by
\begin{equation}\label{eq:ideal-excess}
 \delta_P(\mathbf a)=\max_{I\in J(P)}
 \left(\sum_{e\in I}a_e-|I|\right),
\end{equation}
where $J(P)$ denotes the set of all order ideals of $P$, including the empty ideal. Then $\delta_P(\mathbf a)\ge0$, and
\begin{equation}\label{eq:excess-feasible}
 \mathbf a\in\mathcal Q_P\cap\ZZ^P
 \quad\Longleftrightarrow\quad \delta_P(\mathbf a)=0.
\end{equation}
Write $r_+=\max(r,0)$, and let $\widehat P$ be obtained by adjoining one new minimum below all of $P$.

\begin{lemma}\label{lem:excess}
For a disjoint union $P\sqcup Q$,
\[
 \delta_{P\sqcup Q}(\mathbf a,\mathbf b)
 =\delta_P(\mathbf a)+\delta_Q(\mathbf b).
\]
If the coordinate at the new minimum of $\widehat P$ is $p$, then
\[
 \delta_{\widehat P}(p,\mathbf a)
 =\bigl(p-1+\delta_P(\mathbf a)\bigr)_+.
\]
\end{lemma}

\begin{proof}
An ideal of $P\sqcup Q$ is the union of an ideal of $P$ and an ideal of $Q$, so the maximum separates. A nonempty ideal of $\widehat P$ consists of the new minimum together with an ideal of $P$. Including the empty ideal gives the positive part in the second identity.
\end{proof}

We use Lemma~\ref{lem:excess} for combs in Section~\ref{sec:combs}. For two-row posets, prefix sums give a convenient state description.

Let $a_i,b_j$ be the nonnegative integer coordinates at $e_i,f_j$, respectively. Put
\[
 A_i=\sum_{\ell=1}^i a_\ell,\qquad
 B_j=\sum_{\ell=1}^j b_\ell,\qquad A_0=B_0=0.
\]
An order ideal consists of an initial segment of length $i$ in the first row and one of length $j$ in the second, where $0\le j\le\min(i,k)$. Thus the lattice-point conditions are exactly
\begin{equation}\label{eq:prefix}
 A_i+B_j\le i+j\qquad
 (0\le i\le m,\ 0\le j\le\min(i,k)).
\end{equation}

Process $a_1,b_1,a_2,b_2,\ldots,a_k,b_k$, followed by the remaining first-row coordinates. After some coordinates have been processed, let $i,j$ be their row lengths and set
\[
 p=i-A_i,\qquad q=j-B_j,\qquad
 \mu=\min_{0\le\ell\le j}(\ell-B_\ell),\qquad
 r=p+\mu,\quad s=q-\mu.
\]
The initial state is $(0,0)$. At a first-row step, the second-row prefixes are unchanged, and the inequalities for the new first-row prefix reduce to $r\ge0$. Just before $b_j$ is added, the first row has length $j$ and the minimum $\mu$ is taken over second-row prefixes of lengths $0,\ldots,j-1$. Set $q'=j-B_j$ after this addition. The only new inequality involving an already processed first-row prefix is $p+q'\ge0$, since an ideal containing $f_j$ must contain at least $j$ first-row elements. Together with $p+\mu\ge0$, this is equivalent to $p+\min(\mu,q')\ge0$. Inequalities for longer first-row prefixes are checked when those prefixes are processed. The transitions below therefore enforce all of~\eqref{eq:prefix}.

Adding a first-row coordinate $a$ sends
\begin{equation}\label{eq:lower-step}
 (r,s)\longmapsto(r+1-a,s),\qquad 0\le a\le r+1.
\end{equation}
For a second-row coordinate $b$, put $S=s+1-b$. If $S\ge0$, the state becomes $(r,S)$. If $S<0$, the running minimum decreases by $-S$ and the state becomes $(r+S,0)$; this move is allowed precisely when $r+S\ge0$. Thus the ordinary transitions
\begin{equation}\label{eq:upper-step}
 (r,s)\longmapsto(r,s+1-b),\qquad 0\le b\le s+1,
\end{equation}
are supplemented by
\begin{equation}\label{eq:spill-step}
 (r,s)\longmapsto(r-c,0),\qquad 1\le c\le r,
\end{equation}
corresponding to $b=s+1+c$. Each positive coordinate receives its support weight.

To record the two supports separately, give first-row positive coordinates weight $x$ and second-row positive coordinates weight $y$. On polynomials in $X$, define
\[
 \mathcal C_x(X^r)=X^{r+1}+x\sum_{j=0}^r X^j,
 \qquad
 \mathcal Df(X)=\frac{f(X)-f(1)}{X-1}.
\]
Let $\mathcal E_YF(X,Y)=F(X,1)$. Subscripts indicate the variable on which an operator acts; coefficients in the other variable are left fixed. The first-row and second-row operators on state monomials are, respectively,
\begin{equation}\label{eq:raw-operators}
 \mathcal A=\mathcal C_{x,X},\qquad
 \mathcal B=\mathcal C_{y,Y}+y\mathcal D_X\mathcal E_Y.
\end{equation}
The term $y\mathcal D_X\mathcal E_Y$ records~\eqref{eq:spill-step}, since $\mathcal D_X(X^r)=\sum_{j=0}^{r-1}X^j$. Each transition determines the coordinate that produced it. Operator histories therefore correspond bijectively to feasible coordinate vectors, with the same support weights. Evaluation at $X=Y=1$ sums over all final states.

Make the change of variables $X=1+u$, $Y=1+v$ in $\mathbb Q[x,y][X,Y]$. On $\mathbb Q[x,y][u]$, let $U$ be multiplication by $u$ and let $L$ be the backward shift
\[
 L1=0,\qquad Lu^r=u^{r-1}\quad(r\ge1).
\]
Let $I$ denote the identity operator, let $\mathcal E$ denote evaluation at $1$ in the original variable, and let $E_0$ be projection onto the constant term. Then
\begin{equation}\label{eq:shift-conjugacy}
 \mathcal C_x\longmapsto J_x:=U+(1+x)I+xL,
 \qquad \mathcal D\longmapsto L,
 \qquad \mathcal E\longmapsto E_0.
\end{equation}
For the first identity, use
\[
 (\mathcal C_xf)(X)=Xf(X)+x\frac{Xf(X)-f(1)}{X-1}.
\]
All these operators map polynomials to polynomials, so the change of variables is an algebraic conjugation.

We write $|\vac\rangle=1\otimes1$ and let $\langle\vac|$ extract its coefficient. Operators act from right to left. Since evaluation at $X=Y=1$ becomes evaluation at $u=v=0$, the bivariate support enumerator of $P_{n,n-1}$ is
\begin{equation}\label{eq:biv-moment}
 \mathcal G_n(x,y)=\langle\vac|\bA(\bB\bA)^{n-1}|\vac\rangle,
 \quad
 \bA=J_x\otimes I,\quad
 \bB=I\otimes J_y+yL\otimes E_0.
\end{equation}
In particular, $\mathcal G_n(z,z)=G_n(z)$. For equal weights and general row lengths,
\begin{equation}\label{eq:ferrers-moment}
 H_{P_{m,k}}(z)
 =\langle\vac|\bA^{m-k}(\bB\bA)^k|\vac\rangle,
\end{equation}
where $\bA=J_z\otimes I$ and $\bB=I\otimes J_z+zL\otimes E_0$.

\section{A rank-one transfer identity}\label{sec:trace}

For indeterminates $d,\beta$ and integers $n,r\ge0$, define
\begin{equation}\label{eq:Cnk}
 C_{n,r}(d,\beta)=[w^r](w+d+\beta/w)^n\qquad(r\ge0),
\end{equation}
with $C_{n,r}=0$ for $r>n$. We work over the polynomial ring over $\mathbb Q$ in the displayed parameters. The operators $U,L,E_0$ are as above, and $I$ is the identity. In their monomial matrices, rows index outputs and columns index inputs. We use $[X,Y]=XY-YX$ and write $\CT_v$ for the coefficient of $v^0$. Traces are taken only of matrices with finite support. All power series below are formal.

\begin{theorem}\label{thm:transfer}
Let $d_1,\beta_1,d_2,\beta_2,\kappa$ be independent parameters. Set
\[
 J_i=U+d_iI+\beta_iL,\qquad
 \bA=J_1\otimes I,\qquad
 \bB=I\otimes J_2+\kappa L\otimes E_0,
 \qquad \bT=\bB\bA.
\]
Then, for $n\ge1$,
\begin{equation}\label{eq:general-transfer}
 \langle\vac|\bA\bT^{n-1}|\vac\rangle
 =\frac1n\sum_{r=1}^n\kappa^{r-1}C_{n,r}(d_2,\beta_2)
 \bigl(C_{n,r-1}(d_1,\beta_1)-\beta_1C_{n,r+1}(d_1,\beta_1)\bigr).
\end{equation}
In particular, if $d_1=d_2=d$ and $\beta_1=\beta_2=\kappa=\beta$, then
\begin{equation}\label{eq:transfer-product}
 \langle\vac|\bA\bT^{n-1}|\vac\rangle
 =\frac1nC_{n,1}(d,\beta)C_{n,0}(d,\beta).
\end{equation}
\end{theorem}

\begin{proof}
Let $h_n=\langle\vac|\bA\bT^{n-1}|\vac\rangle$ and $\Pi=E_0\otimes E_0$. Since $[L,U]=E_0$, we have
\[
 [\bB,\bA]=\kappa\Pi,
 \qquad
 [\bA,\bT^n]
 =-\kappa\sum_{j=0}^{n-1}\bT^j\Pi\bA\bT^{n-1-j}.
\]
The summands have finite matrix support. To check this, expand
\[
 \bT=J_1\otimes J_2+
 \kappa(I+d_1L+\beta_1L^2)\otimes E_0.
\]
A transfer changes the first height by an integer between $-2$ and $1$ and the second by at most one. The nonzero matrix entries of $\bT^j\Pi\bA\bT^{n-1-j}$ consequently have output indices at most $(j,j)$ and input indices at most $(2(n-1-j)+1,n-1-j)$, componentwise. The trace of this finite-support rank-one matrix is
\[
 \langle\vac|\bA\bT^{n-1-j}\bT^j|\vac\rangle=h_n.
\]
Taking traces therefore gives
\begin{equation}\label{eq:comm-trace}
 \Tr[\bA,\bT^n]=-\kappa n h_n.
\end{equation}

We compute the same trace by a finite truncation in the first coordinate. Put
\[
 \Delta_n=\bT^n-(J_1\otimes J_2)^n,
 \qquad q_i(v)=v+d_i+\beta_i/v.
\]
The matrix $\Delta_n$ is banded in the first coordinate. Every term in its expansion contains a second-coordinate projection, so its second-coordinate input and output indices are at most $n-1$. Since $\bA$ commutes with $J_1\otimes J_2$, we have $[\bA,\bT^n]=[\bA,\Delta_n]$. Let $(\Delta_n)_{i,j}$ denote its first-coordinate blocks, with blocks having a negative index interpreted as zero. Then
\begin{align*}
 [\bA,\Delta_n]_{i,i}
 ={}&(\Delta_n)_{i-1,i}+\beta_1(\Delta_n)_{i+1,i}\\
 &-(\Delta_n)_{i,i+1}-\beta_1(\Delta_n)_{i,i-1}.
\end{align*}
Summing over $0\le i\le M$ cancels all interior terms and gives
\begin{equation}\label{eq:boundary-blocks}
 \Tr[\bA,\bT^n]
 =\Tr_2\bigl(\beta_1(\Delta_n)_{M+1,M}-(\Delta_n)_{M,M+1}\bigr)
 \qquad(M>2n).
\end{equation}
Here $\Tr_2$ is the trace in the second coordinate; the blocks have finite support in that coordinate. For $M>2n$, the truncated sum includes every nonzero diagonal entry of the commutator.

Starting at first height $M$ or $M+1$, an $n$-transfer path can descend by at most $2n$. If $M>2n$, it stays above the boundary, including at the intermediate stages of each transfer. We may thus replace $U,L$ in the first coordinate by the bilateral shifts $v,v^{-1}$. For the two boundary blocks in~\eqref{eq:boundary-blocks}, this replacement gives the Laurent-polynomial matrix
\[
 q_1(v)^n\bigl((J_2+(\kappa/v)E_0)^n-J_2^n\bigr).
\]
Set
\[
 R_n(c)=\Tr\bigl((J_2+cE_0)^n-J_2^n\bigr).
\]
The difference in this definition has finite matrix support. For this bilateral operator, the $(i,j)$ block is the coefficient of $v^{i-j}$, so the two boundary blocks in~\eqref{eq:boundary-blocks} give
\begin{equation}\label{eq:flux}
 \Tr[\bA,\bT^n]
 =\CT_v(\beta_1/v-v)q_1(v)^nR_n(\kappa/v).
\end{equation}
Only finite-support differences and blocks enter these traces.

It remains to evaluate $R_n$. Define
\[
 f(s)=\langle0|(I-sJ_2)^{-1}|0\rangle,\qquad Z(s)=sf(s).
\]
The first-return decomposition of Motzkin excursions gives
\begin{equation}\label{eq:lagrange-Z}
 f=1+d_2sf+\beta_2s^2f^2,
 \qquad Z=s(1+d_2Z+\beta_2Z^2).
\end{equation}
Put $M_c=J_2+cE_0$. Differentiating the finite-support difference defining $R_n(c)$ gives
\[
 \frac{\partial R_n(c)}{\partial c}
 =\sum_{i=0}^{n-1}\Tr\bigl(M_c^iE_0M_c^{n-1-i}\bigr)
 =n\langle0|M_c^{n-1}|0\rangle.
\]
Each summand has finite matrix support: the factors surrounding $E_0$ are banded. The equality follows from the finite rank-one identity $\Tr(uv^{\mathsf T})=v^{\mathsf T}u$. Defining $\Psi(s,c)=\sum_{n\ge1}R_n(c)s^n/n$, we obtain
\begin{align*}
 \frac{\partial \Psi(s,c)}{\partial c}
 &=s\langle0|(I-s(J_2+cE_0))^{-1}|0\rangle\\
 &=\frac{Z(s)}{1-cZ(s)}.
\end{align*}
The rank-one resolvent identity is
\[
 \langle0|(I-s(J_2+cE_0))^{-1}|0\rangle
 =\frac{f(s)}{1-sc f(s)}.
\]
It gives the second equality above. Integrating coefficientwise in $c$ and using $\Psi(s,0)=0$ yields
\[
 \sum_{n\ge1}\frac{R_n(c)}n s^n=-\log(1-cZ(s)).
\]
Lagrange inversion~\cite[Corollary 5.4.3]{Stanley} applied to~\eqref{eq:lagrange-Z} now gives
\begin{align}
 R_n(c)
 &=[w^{n-1}]\frac{c}{1-cw}(1+d_2w+\beta_2w^2)^n\notag\\
 &=\sum_{r=1}^n c^r C_{n,r}(d_2,\beta_2).\label{eq:rankone-moment}
\end{align}
For the last equality, use
\[
 [w^{n-r}](1+d_2w+\beta_2w^2)^n=C_{n,r}(d_2,\beta_2).
\]
In~\eqref{eq:flux}, the term $(\kappa/v)^r$ contributes
\[
 \kappa^r\bigl(\beta_1C_{n,r+1}(d_1,\beta_1)-C_{n,r-1}(d_1,\beta_1)\bigr).
\]
Comparison with~\eqref{eq:comm-trace} gives~\eqref{eq:general-transfer} after division by $\kappa$ in the fraction field. Both sides of the resulting identity are polynomials in the parameters, so it also holds at $\kappa=0$.

For matching parameters write $C_r=C_{n,r}(d,\beta)$. The sum on the right of~\eqref{eq:general-transfer}, multiplied by $n$, becomes
\[
 \sum_{r=1}^n \beta^{r-1}C_rC_{r-1}
 -\sum_{r=1}^n \beta^rC_rC_{r+1}=C_1C_0.
\]
The intermediate terms cancel, and $C_{n+1}=C_{n,n+1}(d,\beta)=0$, proving~\eqref{eq:transfer-product}.
\end{proof}

By~\eqref{eq:lagrange-Z}, $C_{n,1}(d,\beta)/n$ is the return moment $\langle0|J^{n-1}|0\rangle$ for $J=U+dI+\beta L$. Thus~\eqref{eq:transfer-product} factors into the weighted enumerators of a Motzkin excursion of length $n-1$ and an unrestricted Motzkin bridge of length $n$~\cite{Flajolet}. Horizontal steps, down-steps, and up-steps have weights $d,\beta,1$, respectively.

\section{The Narayana product and its row refinement}\label{sec:refinement}

\begin{proof}[Proof of Theorem~\ref{thm:product}]
Take $d=1+z$, $\beta=z$ in~\eqref{eq:transfer-product} and use~\eqref{eq:biv-moment}. Since
\[
 w+1+z+z/w=(1+w)(1+z/w),
\]
we have
\[
 C_{n,0}(1+z,z)=B_n(z),\qquad
 \frac1nC_{n,1}(1+z,z)=N_n(z).
\]
This proves~\eqref{eq:product}. At $z=1$, Vandermonde's identity gives $B_n(1)=\binom{2n}{n}$ and $N_n(1)=\binom{2n}{n}/(n+1)$, proving~\eqref{eq:total}.
\end{proof}

For $n\ge1$ and $0\le r\le n$, put
\begin{equation}\label{eq:SrDr}
 \begin{split}
 S_{n,r}(z)&=\sum_{j=0}^{n-r}\binom nj\binom n{j+r}z^j,\\
 D_{n,r}(z)&=S_{n,r}(z)-zS_{n,r+2}(z),
 \end{split}
\end{equation}
where $S_{n,r}=0$ for $r>n$. A coefficient comparison shows that
\begin{equation}\label{eq:D-positive}
 D_{n,r}(z)=\frac{r+1}{n+1}
 \sum_{j=0}^{n-r}\binom{n+1}{j}\binom{n+1}{j+r+1}z^j.
\end{equation}

\begin{theorem}\label{thm:bivariate}
For $n\ge1$, let
\[
 \mathcal G_n(x,y)=\sum_{(\mathbf a,\mathbf b)\in \mathcal{Q}_{P_n}\cap\ZZ^{2n-1}}
 x^{|\supp(\mathbf a)|}y^{|\supp(\mathbf b)|},
\]
where the first row has length $n$ and the second has length $n-1$. Then
\begin{equation}\label{eq:bivariate-sum}
 \mathcal G_n(x,y)=\frac1n\sum_{r=0}^{n-1}
 y^r S_{n,r+1}(y)D_{n,r}(x).
\end{equation}
Every summand in~\eqref{eq:bivariate-sum} has nonnegative coefficients.
\end{theorem}

\begin{proof}
Use~\eqref{eq:biv-moment} and Theorem~\ref{thm:transfer} with
\[
 (d_1,\beta_1,d_2,\beta_2,\kappa)=(1+x,x,1+y,y,y).
\]
Then $C_{n,r}(1+z,z)=S_{n,r}(z)$, so~\eqref{eq:general-transfer} is exactly~\eqref{eq:bivariate-sum}. Formula~\eqref{eq:D-positive} follows from
\[
 \binom nj\binom n{j+r}
 -\binom n{j-1}\binom n{j+r+1}
 =\frac{r+1}{n+1}\binom{n+1}{j}\binom{n+1}{j+r+1},
\]
with out-of-range binomial coefficients interpreted as zero.
\end{proof}

For example,
\[
 \mathcal G_2(x,y)=1+3x+x^2+y(2+4x+x^2).
\]
Its two boundary specializations are
\begin{equation}\label{eq:boundary-spec}
 \mathcal G_n(x,0)=N_{n+1}(x),\qquad
 \mathcal G_n(0,y)=\frac1n\sum_{j=0}^{n-1}\binom n{j+1}\binom{2n}{j}y^j.
\end{equation}
For the first identity, only $r=0$ remains in~\eqref{eq:bivariate-sum}; since $S_{n,1}(0)=n$, we obtain $\mathcal G_n(x,0)=D_{n,0}(x)=N_{n+1}(x)$. For the second, use $D_{n,r}(0)=S_{n,r}(0)=\binom nr$. The coefficient of $y^j$ is
\[
 \frac1n\binom n{j+1}\sum_{r=0}^{j}\binom nr\binom n{j-r}
 =\frac1n\binom n{j+1}\binom{2n}{j}.
\]
Setting $y=1$ in the second formula and applying Vandermonde's identity gives $\binom{3n}{n}/(2n+1)$.

\section{Comb posets and an algebraic generating function}\label{sec:combs}

At each step of the comb recursion, call the new minimum a spine element and the other new element a tooth. Write $S_n$ and $L_n$ for the spine elements and teeth of $T_n$, respectively. Each set has $n$ elements. Define
\begin{equation}\label{eq:comb-biv-def}
 \mathcal K_n(a,b)=
 \sum_{\mathbf v\in\mathcal Q_{T_n}\cap\ZZ^{T_n}}
 a^{|\supp(\mathbf v)\cap S_n|}
 b^{|\supp(\mathbf v)\cap L_n|}.
\end{equation}
Thus $\mathcal K_n(z,z)=K_n(z)$ and $\mathcal K_0(a,b)=1$.

\begin{theorem}\label{thm:comb-bivariate}
The series
\[
 W(t;a,b)=t\sum_{n\ge0}\mathcal K_n(a,b)t^n
\]
is the unique solution in $t\mathbb Q[a,b][[t]]$ of
\begin{equation}\label{eq:comb-biv-cubic}
 W=t(1+bW)\bigl(1+(1+b)W\bigr)
       \bigl(1+a(1+b)W\bigr).
\end{equation}
Consequently,
\begin{equation}\label{eq:comb-biv-coeff}
 \mathcal K_n(a,b)=\frac1{n+1}[w^n]
 (1+bw)^{n+1}\bigl(1+(1+b)w\bigr)^{n+1}
 \bigl(1+a(1+b)w\bigr)^{n+1}.
\end{equation}
\end{theorem}

\begin{proof}
For $k\ge0$, define the support enumerator with excess at most $k$ by
\[
 R_{n,k}(a,b)=
 \sum_{\substack{\mathbf v\in\NN^{T_n}\\\delta_{T_n}(\mathbf v)\le k}}
 a^{|\supp(\mathbf v)\cap S_n|}
 b^{|\supp(\mathbf v)\cap L_n|}.
\]
The sum is finite: the constraint for the full ideal $T_n$ gives
$\sum_e v_e\le2n+k$. We have
\begin{equation}\label{eq:comb-slack-initial}
 R_{n,0}=\mathcal K_n,\qquad R_{0,k}=1.
\end{equation}
Let $p$ and $q$ be the new spine and tooth coordinates. By Lemma~\ref{lem:excess},
\begin{equation}\label{eq:comb-excess-step}
 \delta_{T_{n+1}}(p,q,\mathbf v)
 =\bigl(p-1+\delta_{T_n}(\mathbf v)+(q-1)_+\bigr)_+.
\end{equation}
Since $k\ge0$, the condition $\delta_{T_{n+1}}(p,q,\mathbf v)\le k$ is equivalent to
\[
 p+(q-1)_++\delta_{T_n}(\mathbf v)\le k+1.
\]
The contributions of the two new coordinates are recorded by
\[
 \sum_{p\ge0}a^{\mathbf1_{\{p>0\}}}u^p
 =1+\frac{au}{1-u},\qquad
 \sum_{q\ge0}b^{\mathbf1_{\{q>0\}}}u^{(q-1)_+}
 =1+b+\frac{bu}{1-u}.
\]
Write
\begin{equation}\label{eq:comb-kernel-C}
 C(u)=\sum_{j\ge0}c_j(a,b)u^j
 =\frac{(1+(a-1)u)(1+b-u)}{(1-u)^2}.
\end{equation}
Thus $c_0=1+b$ and $c_j=a+b+jab$ for $j\ge1$. Grouping the choices of $(p,q)$ according to $j=p+(q-1)_+$ gives
\begin{equation}\label{eq:comb-slack-rec}
 R_{n+1,k}=
 \sum_{j=0}^{k+1}c_j(a,b)R_{n,k+1-j}.
\end{equation}

Introduce
\[
 \mathcal R(t,u)=\sum_{n,k\ge0}R_{n,k}(a,b)t^nu^k,
 \qquad K(t;a,b)=\mathcal R(t,0)=\sum_{n\ge0}\mathcal K_n(a,b)t^n.
\]
For the kernel calculation, work over $\mathbb Q(a,b)$, with $\mathcal R\in\mathbb Q(a,b)[[t,u]]$. Multiplying~\eqref{eq:comb-slack-rec} by $t^{n+1}u^{k+1}$ and summing over $n,k\ge0$, with the initial values in~\eqref{eq:comb-slack-initial}, gives
\begin{equation}\label{eq:comb-kernel}
 \bigl(u-tC(u)\bigr)\mathcal R(t,u)
 =\frac{u}{1-u}-t(1+b)K(t;a,b).
\end{equation}
There is a unique series $V(t)\in t\mathbb Q(a,b)[[t]]$ satisfying
\begin{equation}\label{eq:comb-small-root}
 V=tC(V).
\end{equation}
This root satisfies $V=(1+b)t+O(t^2)$. Since $V$ has zero constant term, substituting $u=V(t)$ into~\eqref{eq:comb-kernel} involves only finitely many pairs $(n,k)$ in each coefficient of $t$. The kernel vanishes, leaving
\begin{equation}\label{eq:comb-root-relation}
 \frac{V}{1-V}=t(1+b)K(t;a,b)=(1+b)W.
\end{equation}
Hence
\[
 V=\frac{(1+b)W}{1+(1+b)W}.
\]
Substitution into~\eqref{eq:comb-kernel-C} yields
\[
 C(V)=(1+b)(1+bW)\bigl(1+a(1+b)W\bigr).
\]
Equation~\eqref{eq:comb-small-root} gives~\eqref{eq:comb-biv-cubic} after cancellation of $1+b$ over $\mathbb Q(a,b)$. Its right-hand side is $t$ times a polynomial in $W$ with constant term $1$, so it determines $W=t+O(t^2)$ recursively. Lagrange inversion yields~\eqref{eq:comb-biv-coeff}. The kernel substitution was carried out over $\mathbb Q(a,b)$, but the resulting coefficient identities lie in $\mathbb Q[a,b]$ and hence remain valid under specialization.
\end{proof}

For example,
\[
 \mathcal K_1(a,b)=1+a+2b+ab,
\]
which directly counts the five lattice points for the two-element chain $T_1$.

\begin{proof}[Proof of Theorem~\ref{thm:comb}]
Since $F(x,z)=zW(x/z;z,z)$, setting $a=b=z$ and $t=x/z$ in~\eqref{eq:comb-biv-cubic} gives~\eqref{eq:comb-cubic}. Write
\[
 \phi(w,z)=(1+w)\bigl(1+(1+z)w\bigr)
                      \bigl(1+(1+z^{-1})w\bigr).
\]
Lagrange inversion gives
\begin{equation}\label{eq:comb-lagrange}
 K_n(z)=\frac{z^n}{n+1}[w^n]\phi(w,z)^{n+1}.
\end{equation}
The identity
\[
 \phi(w,z)=(1+w)\left(1+\frac{(1+z)^2}{z}w(1+w)\right)
\]
therefore gives
\begin{align*}
 K_n(z)
 &=\frac1{n+1}\sum_{r=0}^n
   \binom{n+1}{r}\binom{n+1+r}{n-r}
   z^{n-r}(1+z)^{2r}\\
 &=\frac1{n+1}\sum_{j=0}^n
   \binom{n+1}{j+1}\binom{2n+1-j}{j}
   z^j(1+z)^{2n-2j},
\end{align*}
where the second equality follows on setting $j=n-r$. This is~\eqref{eq:comb-gamma-main}, and $\eta_{n,j}>0$ for $0\le j\le n$. Integrality follows from $K_n\in\ZZ[z]$: the lowest-degree term of $z^j(1+z)^{2n-2j}$ is $z^j$, with coefficient $1$, so comparison in increasing degree determines each $\eta_{n,j}$ as an integer. The expansion gives the asserted $\gamma$-positivity, symmetry, and unimodality.
\end{proof}

The inverse-series formulation in~\cite[Section 9.6]{AC} follows by setting
\[
 g(w,z)=\sum_{m\ge1}(-1)^{m-1}
 \left(1+(1+z)^m\sum_{j=0}^{m-2}z^j\,z^{-m+1}\right)w^m,
\]
with the inner sum zero when $m=1$. Summing the geometric series gives
\begin{align*}
 g(w,z)
 &=\frac{w}{1+w}
 -\frac{(1+z)(1+z^{-1})w^2}
 {\bigl(1+(1+z)w\bigr)\bigl(1+(1+z^{-1})w\bigr)}\\
 &=\frac{w}{(1+w)\bigl(1+(1+z)w\bigr)
                         \bigl(1+(1+z^{-1})w\bigr)}.
\end{align*}
Thus~\eqref{eq:comb-cubic} is equivalent to $g(F(x,z),z)=x$.

\begin{corollary}\label{cor:comb-total}
For every $n\ge0$,
\begin{equation}\label{eq:comb-total}
 |\mathcal Q_{T_n}\cap\ZZ^{2n}|
 =\frac1{n+1}\sum_{r=0}^n
   \binom{n+1}{r+1}\binom{2n+2}{r}2^r.
\end{equation}
If $A(x)=\sum_{n\ge0}|\mathcal Q_{T_n}\cap\ZZ^{2n}|x^n$, then
\[
 A=(1+xA)(1+2xA)^2.
\]
These lattice-point counts are the small $3$-Schr\"oder numbers, with the indexing used for $A034015$ in~\cite{OEIS}.
\end{corollary}

\begin{proof}
At $z=1$, equation~\eqref{eq:comb-lagrange} becomes
\[
 K_n(1)=\frac1{n+1}[w^n](1+w)^{n+1}(1+2w)^{2n+2}.
\]
Expanding the second factor gives~\eqref{eq:comb-total}. Substituting $F(x,1)=xA(x)$ into~\eqref{eq:comb-cubic} gives the equation for $A$. The sum is a known formula: it is the case $m=2$ of~\cite[Theorem 4.4]{Yan}, with its size parameter replaced by $n+1$. It is also recorded in~\cite{OEIS}; see~\cite[Lemma 3.7]{Easwar}, which derives it from formulas in~\cite{YJ}. Here we identify these numbers with the lattice-point counts for $\mathcal Q_{T_n}$.
\end{proof}

The first three $\gamma$-expansions are
\begin{align*}
 K_1(z)&=(1+z)^2+z,\\
 K_2(z)&=(1+z)^4+4z(1+z)^2+z^2,\\
 K_3(z)&=(1+z)^6+9z(1+z)^4+10z^2(1+z)^2+z^3.
\end{align*}
In particular, $K_n(1)$ starts with $1,5,33,249,2033,17485$ for $n=0,1,\ldots,5$.

\section{Positivity and zeros for two-row posets}\label{sec:consequences}

\begin{proof}[Proof of Theorem~\ref{thm:ferrers-gamma}]
The degree of $H_{P_{m,k}}$ is $m+k$, since the all-ones vector is feasible. For independent indeterminates $d,\beta$, set
\[
 J=U+dI+\beta L,\qquad
 \bA=J\otimes I,\qquad
 \bB=I\otimes J+\beta L\otimes E_0.
\]
These operators specialize to those in~\eqref{eq:ferrers-moment} at $d=1+z$ and $\beta=z$. Their moves on $\NN^2$ are up-steps of weight $1$, horizontal steps of weight $d$, and down-steps of weight $\beta$. The boundary term $\beta L\otimes E_0$ is a down-step in the first coordinate, available when that coordinate is positive and the second is zero.

The operator word $\bA^{m-k}(\bB\bA)^k$ has length $m+k$. Each up-step increases the sum of the two coordinates by one, and each down-step decreases it by one. A return walk with $j$ down-steps therefore has $j$ up-steps and $m+k-2j$ horizontal steps. Its weight is
\[
 \beta^j d^{m+k-2j}.
\]
Thus, for nonnegative integers $\gamma_{m,k,j}$,
\[
 \langle\vac|\bA^{m-k}(\bB\bA)^k|\vac\rangle
 =\sum_{j=0}^{\lfloor(m+k)/2\rfloor}
 \gamma_{m,k,j}\beta^jd^{m+k-2j}.
\]
Substituting $d=1+z$ and $\beta=z$ gives the required $\gamma$-expansion. Here $\gamma_{m,k,j}$ counts return walks with $j$ down-steps, with moves arising from different operator summands counted separately.
\end{proof}

For $P_n$, the $\gamma$-coefficients admit a closed formula. Write $\mathrm{Cat}_j=\binom{2j}{j}/(j+1)$.

\begin{corollary}\label{cor:gamma-explicit}
For $n\ge1$,
\[
 G_n(z)=\sum_{r=0}^{n-1}\gamma_{n,r}z^r(1+z)^{2n-1-2r},
\]
where
\begin{equation}\label{eq:gamma-explicit}
 \gamma_{n,r}=\sum_{\substack{i,j\ge0\\i+j=r}}
 \mathrm{Cat}_i\binom{n-1}{2i}\binom n{2j}\binom{2j}{j}.
\end{equation}
\end{corollary}

\begin{proof}
In the excursion and bridge interpretations following Theorem~\ref{thm:transfer}, choose the positions of the nonhorizontal steps. Deleting the horizontal steps leaves a Dyck path in the excursion case and an unrestricted up-and-down bridge in the other. Hence
\[
 \frac1n C_{n,1}(d,\beta)=\sum_{i=0}^{\lfloor(n-1)/2\rfloor}\mathrm{Cat}_i\binom{n-1}{2i}\beta^id^{n-1-2i},
\]
\[
 C_{n,0}(d,\beta)=\sum_{j=0}^{\lfloor n/2\rfloor}\binom n{2j}\binom{2j}{j}\beta^jd^{n-2j}.
\]
Multiplying these two expansions and setting $d=1+z$ and $\beta=z$ gives~\eqref{eq:gamma-explicit}.
\end{proof}

The product formula also determines the zeros of $G_n$.

\begin{corollary}\label{cor:zeros}
For $n\ge2$, the zeros of $N_n$ strictly interlace those of $B_n$. For every $n\ge1$, $G_n$ has $2n-1$ simple negative roots. Its coefficients are strictly log-concave for $n\ge2$.
\end{corollary}

\begin{proof}
Let $\mathcal L_n$ be the Legendre polynomial and put $\xi=(1+z)/(1-z)$. The standard binomial formula for Jacobi polynomials and their derivative relation~\cite[18.5.8 and 18.9.15]{DLMF} give
\[
 B_n(z)=(1-z)^n\mathcal L_n(\xi),\qquad
 N_n(z)=\frac{2}{n(n+1)}(1-z)^{n-1}\mathcal L_n'(\xi).
\]
The polynomial $\mathcal L_n$ has $n$ simple roots in $(-1,1)$; see~\cite[Section 18.2(vi)]{DLMF}. Its derivative has degree $n-1$, so Rolle's theorem gives exactly one root between each consecutive pair. The two root sets therefore strictly interlace. The inverse transformation
\[
 z=\frac{\xi-1}{\xi+1}
\]
is strictly increasing from $(-1,1)$ onto $(-\infty,0)$ and preserves both simplicity and interlacing. The product $G_n=N_nB_n$ therefore has $2n-1$ simple negative roots. Write $G_n(z)=\sum_{j=0}^{D}a_jz^j$, where $D=2n-1$. Its coefficients are positive, and Newton's inequalities give
\[
 a_j^2\ge\frac{(j+1)(D-j+1)}{j(D-j)}a_{j-1}a_{j+1}
 >a_{j-1}a_{j+1}\qquad(1\le j\le D-1).
\]
This proves strict log-concavity for $n\ge2$.
\end{proof}

For example,
\[
 H_{P_{3,1}}(z)=1+8z+15z^2+8z^3+z^4
 =(1+z)^4+4z(1+z)^2+z^2.
\]
Theorem~\ref{thm:ferrers-gamma} applies both to such unequal row lengths and to the full grids $P_{n,n}=C_n\times C_2$. The product evaluation concerns the subfamily $P_{n,n-1}$.

\section*{Acknowledgments}
This work was supported by the National Natural Science Foundation of China (No. 12401441) and the Natural Science Foundation of Hunan Province (No. 2025JJ60010).

\end{document}